# A NOTE ON MULTITYPE BRANCHING PROCESSES WITH IMMIGRATION IN A RANDOM ENVIRONMENT

By Alexander Roitershtein

*University of British Columbia*

We consider a multitype branching process with immigration in a random environment introduced by Key in [*Ann. Probab.* **15** (1987) 344–353]. It was shown by Key that, under the assumptions made in [*Ann. Probab.* **15** (1987) 344–353], the branching process is subcritical in the sense that it converges to a proper limit law. We complement this result by a strong law of large numbers and a central limit theorem for the partial sums of the process. In addition, we study the asymptotic behavior of oscillations of the branching process, that is, of the random segments between successive times when the extinction occurs and the process starts again with the next wave of the immigration.

**1. Introduction and statement of results.** In this paper we consider the *multitype branching process with immigration in a random environment* (MBPIRE) introduced by Key in [14]. Broadly speaking, the MBPIRE is a vector-valued random process $Z_n = (Z_n^{(1)}, Z_n^{(2)}, \ldots, Z_n^{(d)})$ that describes evolution of a population of particles of $d$ different types in discrete time $n = 0, 1, \ldots$. The integer number $d \geq 1$ is fixed and $Z_n^{(i)}$ denotes the number of particles of type $i$ in generation $n$. We next present the process $Z = (Z_n)_{n \geq 0}$ in more detail.

The *environment* $\omega = (\omega_n)_{n \in \mathbb{Z}}$ is a stationary ergodic sequence of random variables $\omega_n$ taking values in a measurable space $(\mathcal{S}, \mathcal{B})$. We denote by $P$ the law of $\omega$ on $(\Omega, \mathcal{F}) := (\mathcal{S}^{\mathbb{Z}}, \mathcal{B}^{\otimes \mathbb{Z}})$ and by $E_P$ the expectation with respect to $P$. Both the law of the immigration and the branching mechanism depend on the realization of the environment.

Let $\mathbb{Z}_+$ denote $\mathbb{N} \cup \{0\}$. The *immigration* $X = (X_n)_{n \in \mathbb{Z}_+}$ is a sequence of independent conditionally on $\omega$ random vectors $X_n = (X_n^{(1)}, \ldots, X_n^{(d)}) \in \mathbb{Z}_+^d$, describing the number of immigrants of all types arrived at time $n \in \mathbb{Z}_+$. Let $\mathcal{P}_d$ be the set of probability measures on $\mathbb{Z}_+^d$ and assume that the law $q_{\omega,n}$ of $X_n$ in a fixed environment $\omega$ is a function of the past state $\omega_{n-1}$ of the environment:

$$q_{\omega,n} = q(\omega_{n-1}) \tag{1.1}$$

for some measurable function $q : \mathcal{S} \to \mathcal{P}_d$. Here we regard $\mathcal{P}_d$ as a subspace of the space of bounded real-valued sequences $l_\infty$ endowed with the Borel $\sigma$-algebra induced by the product topology.

Every particle in the branching process is a descendant of an immigrant but the immigrants themselves are not counted in the population. At each unit of time, every particle of type $i$ present in the system (*including the immigrants*) splits (for a given realization of the environment independently of the previous history and of the other particles) into a random number of offspring of all types according to a law $p_{\omega,n}^{(i)}$ which is a function of $\omega_n$:

$$p_{\omega,n}^{(i)} : p_i(\omega_n) \tag{1.2}$$

for some measurable functions $p_i : \mathcal{S} \to \mathcal{P}_d, i = 1, \ldots, d$.

Thus, conditionally on the environment, $(Z_n)_{n \in \mathbb{Z}_+}$ is a nonhomogeneous Markov chain that satisfies the initial condition $Z_0 = \mathbf{0}$ [here and throughout $\mathbf{0}$ stands for $(0, \ldots, 0) \in \mathbb{Z}_+^d$] and the branching equation

$$Z_{n+1} = \sum_{i=1}^{d} \sum_{m=1}^{Z_n^{(i)} + X_n^{(i)}} L_{n,m}^{(i)} \qquad \text{for } n \geq 0, \tag{1.3}$$

where the random vectors $L_{n,m}^{(i)} = (L_{n,m}^{(i,1)}, \ldots, L_{n,m}^{(i,d)}) \in \mathbb{Z}_+^d$ (with $n \in \mathbb{Z}_+; i = 1, \ldots, d; m \in \mathbb{N}$) are independent and $L_{n,m}^{(i)}$ are distributed according to the law $p_{\omega,n}^{(i)}$ for every $m \in \mathbb{N}$.

The underlying probability space $(\Omega \times \mathcal{T}, \mathcal{F} \otimes \mathcal{G}, \mathbb{P})$, where $\mathcal{T}$ is a space of family trees that describes the "genealogy" of the particles and $\mathcal{G}$ is its Borel $\sigma$-algebra, is constructed using the recipe given in [12], Chapter VI. The measure $\mathbb{P}$ on $\Omega \times \mathcal{T}$ is defined as $P \otimes P_\omega$, where for $\omega \in \Omega$, $P_\omega$ is the (popularly known as *quenched*) law on $(\mathcal{T}, \mathcal{G})$ consistent with the preceding description of the process in the environment $\omega$. The marginal $\mathbb{P}(\cdot) = E_P(P_\omega(\cdot))$ of $\mathbb{P}$ on $(\mathcal{T}, \mathcal{G})$ is referred to as the *annealed* law of the process. The expectations with respect to $P_\omega$ and $\mathbb{P}$ are denoted by $E_\omega$ and $\mathbb{E}$, respectively. For a random vector $v = (v^{(1)}, \ldots, v^{(d)}) \in \mathbb{R}^d$, $E_P(v)$ [correspondingly $E_\omega(v)$, $\mathbb{E}(v)$] stands for the vector whose $i$th component is $E_P(v^{(i)})(E_\omega(v^{(i)}), \mathbb{E}(v^{(i)}))$.

To formulate our assumptions on the process $(Z_n)_{n \in \mathbb{Z}_+}$ we need to introduce some additional notation. For any constant $\beta > 0$ and real-valued random variable $X$ we set

$$\|X\|_{\omega,\beta} := (E_\omega(|X|^\beta))^{1/\beta}. \tag{1.4}$$

We often use $L_{n,1}^{(i)}$ as "a random vector distributed according to the law $p_{\omega,n}^{(i)}$" and in such cases omit the lower index 1 writing it simply as $L_n^{(i)}$. The $j$th component of the vector $L_n^{(i)}$ is denoted by $L_n^{(i,j)}$. For any $n \in \mathbb{Z}_+$ and $i, j = 1, \ldots, d$ we define the following functions of the environment $\omega$:

$$I_{n,\beta}^{(j)}(\omega) := \|X_n^{(j)}\|_{\omega,\beta}, \qquad M_{n,\beta}^{(i,j)}(\omega) := \|L_n^{(i,j)}\|_{\omega,\beta}.$$

We denote by $I_{n,\beta}$ the vector $(I_{n,\beta}^{(1)}, \ldots, I_{n,\beta}^{(d)})$ and (with a slight abuse of notation) by $M_{n,\beta}$ the matrix whose $i$th *column* (not the row) is the vector $(M_{n,\beta}^{(i,1)}, \ldots, M_{n,\beta}^{(i,d)})$. To simplify the notation we write

$$I_n := I_{n,1} \quad \text{and} \quad M_n := M_{n,1}.$$

Finally, for $v = (v^{(1)}, \ldots, v^{(d)}) \in \mathbb{R}^d$ we set $\|v\|_\beta =: (\sum_{i=1}^d |v^{(i)}|^\beta)^{1/\beta}$, and for a $d \times d$ matrix $A$ denote by $\|A\|_\beta$ the corresponding operator norm.

The following is the basic set of conditions imposed in this paper.

ASSUMPTION 1.1.  *We have:*

(A1) $E_P(\log^+ \|I_0\|_1) < \infty$, *where* $\log^+ x := \max\{0, \log x\}$.
(A2) $E_P(\log^+ \|M_0\|_1) < \infty$.
(A3) $\lim_{n\to\infty} n^{-1} E_P(\log \|M_{n-1} \cdots M_1 M_0\|_1) < 0$.

Note that (A1) implies $P(I_0^{(i)} < \infty) = 1$ and $P(M_0^{(i,j)} < \infty) = 1$ for all $i, j$, and that (A2) yields the existence of the limit in (A3) by the sub-additive ergodic theorem of [11].

The model of MBPIRE was introduced by Key in [14], where it is also shown that the process is subcritical under Assumption 1.1, that is, $Z_n$ converges in distribution to a proper limit law (see Theorem 2.1 below where a more precise statement is cited from [14]).

The model is a generalization on one side for the multitype branching processes in random environment (without immigration) introduced by Athreya and Karlin in [3] and on the other side for a special single-type case considered by Kesten, Kozlov and Spitzer in [13]. The latter process is closely related to one-dimensional random walks in random environment (RWRE) and from this point of view was studied by many authors (see, e.g., [1, 8, 13] and Section 2 of the survey [17]). The limiting distribution of $Z_n$ [$\pi_\omega$ defined below in (4.1) turns out to be a geometric law in this case], the *weak* law of large numbers and the central limit theorem for the partial sums of this process in a stationary ergodic environment can be obtained by using corresponding statements for the associated RWRE.

In this paper we study the asymptotic behavior of the sequence $(Z_n)_{n \in \mathbb{Z}_+}$ exploiting both tools developed in the theory of random motion in random media (e.g., the method of the "environment seen from the position of the particle") as well as the construction of the limiting distribution for the general MBPIRE provided in [14].

The general idea of the proofs is as follows. Let $U_{i,n} \in \mathbb{R}^d$ represent the progeny in generation $n$ of all immigrants arrived at time $i$. Then $Z_n = \sum_{i=0}^{n-1} U_{i,n}$ stochastically grows to its stationary distribution as $n$ goes to infinity. The limiting law corresponds to the random variable $\widetilde{Z}_n = \sum_{i=-\infty}^{n-1} U_{i,n}$ which is a.s. finite under Assumption 1.1. Here we introduce "demo" immigrants arriving at negative times and keep the notation $U_{i,n}$ for their progeny. The process $\widetilde{Z}_n$ is stationary and moreover, it turns out that it inherits some ergodic and mixing properties of the environment. Furthermore, since the MPBREs $(U_{i,n})_{n \geq i}$ extinct a.s. for any $i \in \mathbb{Z}$ (cf. [3]), the random variables $Z_n$ and $\widetilde{Z}_n$ coincide for all $n$ large enough. Therefore one can obtain limit theorems for the process $Z_n$ from standard results for the ergodic (mixing under additional assumptions) stationary process $\widetilde{Z}_n$. Note that in this aspect the process with immigration is quite different from the supercritical MPBRE without immigration, where $Z_n$ grows in exponential rate (cf. [6]).

We next turn to the presentation of our main results. Let

$$S_n = \sum_{i=0}^{n-1} Z_n.$$

The following strong law of large numbers is proved in Section 2:

THEOREM 1.2. *Let Assumption* 1.1 *hold. Then* $\mathbb{P}(\lim_{n \to \infty} \frac{S_n}{n} = \rho) = 1$, *where*

(1.5) $$\rho := \sum_{n=0}^{\infty} E_P(M_n \cdots M_1 M_0 I_0) \in (\mathbb{R}_+ \cup \{\infty\})^d.$$

The proof of the theorem is by coupling $(Z_n)_{n \in \mathbb{Z}_+}$ with a stationary and ergodic process $(\widetilde{Z}_n)_{n \in \mathbb{Z}}$ [introduced below in (2.2)] such that $\widetilde{Z}_n$ are distributed according to the limit law of $Z_n$, and applying then the ergodic theorem to $(\widetilde{Z}_n)_{n \in \mathbb{Z}}$ (cf. Lemma 2.2).

In the case where $(\omega_n)_{n \in \mathbb{Z}}$ is a strongly mixing sequence with exponential rate we obtain in Section 3 the following central limit theorem. Let

(1.6) $\alpha(n) = \sup\{|P(A \cap B) - P(A)P(B)| : A \in \mathcal{F}^n, B \in \mathcal{F}_0\}$,

where, for $n \in \mathbb{Z}$,

(1.7) $\mathcal{F}^n := \sigma(\omega_i : i \geq n)$ and $\mathcal{F}_n := \sigma(\omega_i : i < n).$

Recall that the sequence $\omega_n$ is called *strongly mixing* if $\alpha(n) \to_{n \to \infty} 0$. For strongly mixing sequence it holds that (cf. [9], page 10)

(1.8) $|E_P(fg) - E_P(f)E_P(g)| \leq 4\alpha(n)$

for all $\mathcal{F}^n$-measurable random variables $f(\omega)$ and all $\mathcal{F}_0$-measurable random variables $g(\omega)$ such that $P(|f| \leq 1 \text{ and } |g| \leq 1) = 1$.

For any $n \in \mathbb{Z}_+$ and $i \in \{1, \ldots, d\}$ let

$$Y_n^{(i)} := \#\{\text{progeny of type } i \text{ over all generations}$$
$$\text{of all } X_n \text{ immigrants arrived at time } n\}$$

and denote $Y_n := (Y_n^{(1)}, \ldots, Y_n^{(d)}) \in \mathbb{Z}_+^d$. An equivalent definition of $Y_n$ is given in (3.1) below.

THEOREM 1.3. *Let Assumption* 1.1 *hold and suppose in addition that:*

(B1) $\limsup_{n \to \infty} n^{-1} \log E_P(\|M_{n-1,\kappa} \cdots M_{0,\kappa} I_{0,\kappa}\|_\kappa^\kappa) < 0$ *for some* $\kappa > 2$.

(B2) $\limsup_{n \to \infty} n^{-1} \log \alpha(n) < 0$, *where the mixing coefficients* $\alpha(n)$ *are defined in* (1.6).

*Then* $n^{-1/2}(S_n - n \cdot \rho)$ *converges in distribution to a Gaussian random vector* $S_\infty$ *with zero mean. Moreover, if the following condition holds for some* $i \in \{1, \ldots, d\}$,

(B3$^{(i)}$) *There is no function* $f_i : \Omega \to \mathbb{Z}_+^d$ *such that* $P_\omega(Y_0^{(i)} = f_i(\omega)) = 1$ *for $P$-a.e. environment* $\omega$,

*then* $\mathbb{E}([S_\infty^{(i)}]^2) > 0$ *for that* $i$.

REMARK 1.4. (i) Moment condition (B1) is used in the proof of Theorem 1.3 to ensure that $\mathbb{E}(\|Y_0\|_\kappa^\kappa) < \infty$. Mimicking the estimates on the moments of a single-type branching process carried out in the course of the proof of [8], Lemma 2.4 (see also [17], Lemma 2.4.16), it is not hard to check that the conclusion of Theorem 1.3 still holds if (B1) is replaced by

(B1′) There exist (nonrandom) constants $\kappa > 2$ and $m > 0$ such that $P(\|X_0\|_{\omega,\kappa} < m) = 1$, $P(\|L_0^{(i)}\|_{\omega,\kappa} < m) = 1$ for all $i = 1, \ldots, d$, and in addition $\limsup_{n \to \infty} n^{-1} \log E_P(\|M_{n-1} \cdots M_0\|_\kappa^\kappa) < 0$.

Note that (B1′) implies by Jensen inequality that the conditions of Theorem 1.2 hold.

(ii) The random vectors $Y_n$, $n \in \mathbb{Z}_+$, are independent in a fixed environment, and one can check that if Assumption 1.1 holds and condition (B1) is satisfied, then the Lindeberg condition (cf. [10], page 116) is fulfilled for the sequence $(Y_n \cdot t)_{n \in \mathbb{Z}_+}$ for every $t \neq \mathbf{0} \in \mathbb{R}^d$. Therefore, for $P$-a.e. $\omega$, the sequence $(Y_n \cdot t)_{n \in \mathbb{Z}_+}$ obeys a CLT with a random centering under the quenched measure $P_\omega$. Alili (cf. [1], Theorem 5.1 and Section 6) provided a set of conditions on the environment [excluding i.i.d. sequences $(\omega_n)_{n \in \mathbb{Z}_+}$] sufficient to replace the random centering in the quenched CLT of this kind by a constant.

This theorem is similar to the CLT for hitting times of one-dimensional random walks in a stationary ergodic environment (cf. [17], see also [2], Theorem 4.3) that can be equivalently restated in terms of the single-type process considered in [13]. The analogy between hitting times of RWRE and random variables $Y_n$ that exists in the special case can be carried over, to some extent, to the general MBPIRE. Again, the limiting distribution of $Z_n$ is an important ingredient of the proof.

Our next result deals with the regeneration times $\nu_n$ defined by

$$\nu_0 = 0 \quad \text{and} \quad \nu_n = \inf\{i > \nu_{n-1} : Z_i = \mathbf{0}\}, \tag{1.9}$$

with the usual convention that the infimum over an empty set is $\infty$.

THEOREM 1.5. *Let Assumption* 1.1 *hold and suppose in addition that there exists constant $m \in \mathbb{N}$ such that*

$$P(P_\omega(Z_m = \mathbf{0} | Z_0 = \mathbf{1}) > 0) > 0, \tag{1.10}$$

*where $\mathbf{1} := (1, \ldots, 1) \in \mathbb{R}^d$.*

*Then $\mathbb{P}(Z_n = \mathbf{0} \text{ i.o.}) = 1$. Moreover, if*

$$P(P_\omega(Z_1 = \mathbf{0} | Z_0 = \mathbf{1}) > 0) = 1, \tag{1.11}$$

*then $\mathbb{P}(\lim_{n \to \infty} \nu_n / n = \mu) = 1$ for some constant $\mu > 0$.*

The proof of this theorem given in Section 4.1 uses the same coupling with the sequence $\widetilde{Z}_n$ as in the proof of Theorem 1.2 along with a change of measure argument. Assumption (1.10) is essentially a condition of Key implying that the limiting distribution of $Z_n$ puts a positive probability on $\mathbf{0}$ (see Lemma 4.2 below).

The law of large numbers for $\nu_n$ is derived from a somewhat more general result stated in Proposition 4.3. Let $Z_{[i,j]}$ denote the segment of the branching process between times $i$ and $j$, that is,

$$Z_{[i,j]} := (Z_i, Z_{i+1}, \ldots, Z_j), \qquad i \leq j, \tag{1.12}$$

and let $\theta$ denote the shift operator on the space of environments $\Omega$, that is,

$$(\theta\omega)_n = \omega_{n+1}. \tag{1.13}$$

We introduce in Section 4 a probability measure $\widetilde{\mathbb{P}}$ on the underlying measurable space $(\Omega \times \mathcal{T}, \mathcal{F} \times \mathcal{G})$ such that under $\widetilde{\mathbb{P}}$ the sequence of triples $((\theta^{\nu_{n-1}}\omega, \nu_n - \nu_{n-1}, Z_{[\nu_{n-1}+1, \nu_n]}))_{n \in \mathbb{N}}$ is stationary and ergodic. The distribution of the sequence $(Z_n)_{n \geq 0}$ under $\widetilde{\mathbb{P}}$ is the Palm measure $\widetilde{\mathbb{P}}((Z_n)_{n \geq 0} \in \cdot) = \mathbb{P}((\widetilde{Z}_n)_{n \geq 0} \in \cdot | \widetilde{Z}_0 = \mathbf{0})$ and it turns out to be equivalent to $\mathbb{P}((Z_n)_{n \geq 0} \in \cdot)$ if (1.11) holds.

We next consider the case where the law of the branching process is "uniformly elliptic" and the environment $\omega$ is a uniformly mixing sequence. In this case we show that the sequence of pairs $((\nu_n - \nu_{n-1}, Z_{[\nu_{n-1}+1,\nu_n]}))_{n\in\mathbb{N}}$ is uniformly mixing and converges in law to its stationary distribution. Let

$$\varphi(n) = \sup\{|P(A|B) - P(A)| : A \in \mathcal{F}^n, B \in \mathcal{F}_0, P(B) > 0\},$$

where the $\sigma$-algebras $\mathcal{F}^n$ and $\mathcal{F}_0$ are defined in (1.7). Recall that the stationary sequence $(\omega_n)_{n\in\mathbb{Z}}$ is called *uniformly mixing* if $\varphi(n) \to_{n\to\infty} 0$. For uniformly mixing sequence it holds that (cf. [9], page 9)

(1.14) $$|E_P(fg) - E_P(f)E_P(g)| \leq 2\varphi(n)E_P(|f|)$$

for all $\mathcal{F}^n$-measurable random variables $f(\omega)$ and all $\mathcal{F}_0$-measurable random variables $g(\omega)$ such that $P(|f| \leq 1 \text{ and } |g| \leq 1) = 1$ [cf. with (1.8)].

The following theorem is proved in Section 4.2.

THEOREM 1.6. *Let Assumption* 1.1 *hold and suppose in addition that the sequence* $(\omega)_{n\in\mathbb{Z}}$ *is uniformly mixing and that there exist* $d+1$ *measures* $e^{(i)} \in \mathcal{P}_d$, $i = 1, \ldots, d+1$, *and constant* $\varepsilon \in (0,1)$ *such that the following hold:*

(i) *For every* $v \in \mathbb{Z}_+^d$,

$$P(p_{\omega,0}^{(i)}(v) > \varepsilon e^{(i)}(v) \text{ for } i = 1,\ldots,d \text{ and } q_{\omega,0}(v) > \varepsilon e^{(d+1)}(v)) = 1,$$

(ii) $e^{(i)}(\mathbf{0}) > 0$ *for all* $i = 1, \ldots, d$.

*Then:*

(a) *The sequence* $((\nu_n - \nu_{n-1}, Z_{[\nu_{n-1}+1,\nu_n]}))_{n\in\mathbb{N}}$ *is uniformly mixing. More precisely,*

$$\sup_{m\in\mathbb{N}} \sup_{A\in\mathcal{G}^{\nu_n+m}} \sup_{B\in\mathcal{G}_{\nu_m},\mathbb{P}(B)>0} \{|\mathbb{P}(A|B) - \mathbb{P}(A)|\} \to_{n\to\infty} 0,$$

*where we denote* $\mathcal{G}^n = \sigma(Z_i : i \geq n)$ *and* $\mathcal{G}_n = \sigma(Z_i, \omega_i : i < n)$.

(b) *The sequence* $((\nu_n - \nu_{n-1}, Z_{[\nu_{n-1}+1,\nu_n]}))_{n\in\mathbb{N}}$ *converges to its stationary distribution.*

Note that the assumptions of the theorem imply that condition (1.11) of Theorem 1.5 is valid.

The rest of the paper is organized as follows. Section 2 contains the proof of the law of large numbers (Theorem 1.2), Section 3 is devoted to the proof of the central limit theorem (Theorem 1.3), and the proofs of Theorems 1.5 and 1.6 are included in Section 4.

**2. The LLN: Proof of Theorem 1.2.** Theorem 1.2 is obtained in this section from the corresponding statement for the stationary and ergodic "version" $(\widetilde{Z})_{n \in \mathbb{Z}}$ of the process $(Z_n)_{n \in \mathbb{Z}_+}$ which is introduced below in (2.2).

For any $i \in \mathbb{Z}_+$ and $n \geq i$ define the vector $U_{i,n} = (U_{i,n}^{(1)}, \ldots, U_{i,n}^{(d)})$ by setting $U_{i,i} = X_i$ and, for $n > i$,

(2.1)
$$U_{i,n}^{(j)} = \#\{\text{progeny of type } j \text{ in generation } n$$
$$\text{of all immigrants arrived at time } i\}.$$

Thus $Z_n = \sum_{i=0}^{n-1} U_{i,n}$ for $n \in \mathbb{Z}_+$. The sequence $\mathcal{U}_i := (U_{i,n})_{n \geq i}$ form a multitype branching process in random environment as introduced by Athreya and Karlin in [3] (abbreviated as MBPRE in what follows). Assume that the underlying probability space is enlarged to include random vectors $X_i$ and processes $\mathcal{U}_i$ for $i < 0$ such that:

(i) For every $i \in \mathbb{Z}$, $U_{ii} = X_i$.
(ii) The processes $(\mathcal{U}_i)_{i \in \mathbb{Z}}$ are independent in a fixed environment.
(iii) For any $i < 0$, $P_\omega(\mathcal{U}_i \in \cdot) = P_{\theta^i \omega}(\mathcal{U}_0 \in \cdot)$, where the shift $\theta$ is defined in (1.13).

That is, $(\mathcal{U}_i)_{i \in \mathbb{Z}}$ is a stationary collection of MBPREs with independent lines of descent in a fixed environment, each $\mathcal{U}_i$ is starting from the immigration wave described by $X_i$.

Let

(2.2)
$$\widetilde{Z}_n = \sum_{i=-\infty}^{n-1} U_{i,n}, \qquad n \in \mathbb{Z}.$$

Some components of the vectors $\widetilde{Z}_n$ may be infinite a priori. However, the following was in fact proved by Key (cf. [14], Theorem 3.3; the random variables $\widetilde{Z}_n$ are not defined explicitly in [14] but they are "present implicitly," e.g., in Lemma 2.2 there).

THEOREM 2.1 ([14]). *Let Assumption 1.1 hold. Then $\widetilde{Z}_0$ is a proper random vector, that is, $\mathbb{P}(\widetilde{Z}_0^{(i)} < \infty \ \forall i = 1, \ldots, d) = 1$, and $\mathbb{P}(\widetilde{Z}_0 = v) = \lim_{n \to \infty} \mathbb{P}(Z_n = v)$ for every $v \in \mathbb{Z}_+^d$.*

That is, $(\widetilde{Z}_n)_{n \in \mathbb{Z}}$ provides a "stationary version" of the process $Z_n$, and the latter converges in law to its stationary distribution. We next show that the sequence $(\widetilde{Z}_n)_{n \in \mathbb{Z}}$ is ergodic and derive the law of large numbers for the partial sums of $Z_n$ from this fact.

LEMMA 2.2. *Let Assumption 1.1 hold. Then:*

(a) *The sequence $\mathbb{U} = (\mathcal{U}_n)_{n \in \mathbb{Z}}$ is ergodic.*
(b) *$(\widetilde{Z}_n)_{n \in \mathbb{Z}}$ is a stationary ergodic sequence.*
(c) *$\mathbb{P}(\exists\, n_0 \in \mathbb{N} : \widetilde{Z}_n = Z_n \text{ for } n \geq n_0) = 1$.*
(d) *For any measurable function $f : \mathbb{Z}_+^d \to \mathbb{R}$ the following holds $\mathbb{P}$-a.s., provided that the expectation in the right-hand side exists: $\frac{1}{n} \sum_{i=0}^{n-1} f(Z_i) \to_{n \to \infty} \mathbb{E}(f(\widetilde{Z}_0))$.*

PROOF. (a) Let $A \in \sigma(\mathcal{U}_n : n \in \mathbb{Z})$ be an invariant set of the shift operator $\Theta$ defined by $(\Theta \mathbb{U})_n = \mathcal{U}_{n+1}$, that is, $A = \Theta^{-1} A$ modulo a $\mathbb{P}$-null set. Then $P$-a.s., $P_\omega(A) = P_\omega(\Theta^{-1} A) = P_{\theta^{-1}\omega}(A)$. Therefore the function $f(\omega) := P_\omega(A)$ is invariant under the ergodic shift $\theta$ on the sequence $\omega = (\omega_n)_{n \in \mathbb{Z}}$, that is, $f(\omega) = f(\theta\omega)$, $P$-a.s. Hence $f(\omega)$ is a $P$-a.s. constant function. Since the random variables $\mathcal{U}_n$ are independent under $P_\omega$, Kolmogorov's 0–1 law implies that $f(\omega) \in \{0, 1\}$. It follows that $\mathbb{P}(A) = E_P(f(\omega)) \in \{0, 1\}$.

(a) The claim follows from part (a) of the lemma and the definition (2.2) of $\widetilde{Z}_n$.

(c) Follows from the extinction criterion for MBPRE given in [3], Theorem 12 (with notation of Remark following this theorem). The criterion is applied to the progeny of all $\widetilde{Z}_0$ particles living at time 0 in the process $(\widetilde{Z}_n)_{n \in \mathbb{Z}}$, and implies that the extinction of the branching process formed by these particles occurs eventually under Assumption 1.1.

(d) Follows from the ergodic theorem applied to the sequence $(\widetilde{Z}_n)_{n \in \mathbb{Z}}$ and part (c) of the lemma. □

Part (d) of the lemma contains the claim of Theorem 1.2 and we now turn to the proof of the CLT for the partial sums of $Z_n$.

**3. The CLT: Proof of Theorem 1.3.** The CLT for the process $Z_n$ is derived here from the corresponding statement for the sequence $Y_n$ (introduced right before the statement of Theorem 1.6). The definition of $Y_n$ can be equivalently written as follows:

$$(3.1) \qquad Y_n = \sum_{i=n+1}^{\infty} U_{n,i}, \qquad n \in \mathbb{Z},$$

where the random vectors $U_{n,i}$ are introduced in (2.1). It is not hard to check that each component of the vector $E_\omega(Y_0) = \sum_{i=1}^{\infty} M_{i-1} \cdots M_0 I_0$ is $P$-a.s. finite under Assumption 1.1 (see, e.g., [14], Lemma 3.1). Therefore $Y_0$ is a proper random vector, $\mathbb{E}(Y_0) = \rho$, and it follows from Lemma 2.2(a) that $(Y_n)_{n \in \mathbb{Z}}$ is a stationary ergodic sequence. In order to prove the CLT for the partial sums of $Y_n$ we shall use the following general CLT for mixing vector-valued sequences. Recall the definition of the norm $\|v\|_\kappa = (\sum_{i=1}^d |v^{(i)}|^\kappa)^{1/\kappa}$, $v \in \mathbb{R}^d$.

LEMMA 3.1. *Let $(Y_n)_{n\in\mathbb{Z}_+}$ be a stationary sequence of random vectors $Y_n = (Y_n^{(1)}, \ldots, Y_n^{(d)})$ in $\mathbb{R}^d$ such that $\mathbb{E}(\|Y_0\|_\kappa^\kappa) < \infty$ for some $\kappa > 2$. For $n \in \mathbb{N}$ let $\mathcal{F}_n = \sigma(Y_m : m \geq n)$ and denote*

(3.2) $\quad \chi_n = \sup\{|\mathbb{P}(A \cap B) - \mathbb{P}(A)\mathbb{P}(B)| : A \in \sigma(Y_0), B \in \mathcal{F}_n\}.$

*Let $\overline{Y}_n = Y_n - \mathbb{E}(Y_0)$ and $\overline{T}_n = \sum_{i=0}^{n-1} \overline{Y}_i$.*

*If $\sum_{n=1}^\infty \chi_n^{(\kappa-2)/(2\kappa)} < \infty$, then $\overline{T}_n/\sqrt{n}$ converges in distribution to a Gaussian random vector $T_\infty = (T_\infty^{(1)}, \ldots, T_\infty^{(d)})$ with zero mean such that*

(3.3) $\quad \mathbb{E}(T_\infty^{(i)} T_\infty^{(j)}) = \mathbb{E}(\overline{Y}_0^{(i)} \overline{Y}_0^{(j)}) + \sum_{n=1}^\infty \mathbb{E}(\overline{Y}_0^{(i)} \overline{Y}_n^{(j)} + \overline{Y}_0^{(j)} \overline{Y}_n^{(i)}),$

*where the last series converges absolutely.*

The lemma is a combination of a one-dimensional CLT for mixing sequences (cf. [10], page 425) with the multidimensional Cramér–Wold device (cf. [10], page 170).

The next lemma is similar in spirit to [17], Lemma 2.1.10 and the proof of the latter works nearly verbatim.

LEMMA 3.2. *Let the conditions of Theorem 1.3 hold. Then*

$$\limsup_{n\to\infty} 1/n \log \chi_n < 0,$$

*where $\chi_n$ are defined by (3.2).*

PROOF. Let $\mathfrak{Y}_n$ denote the sequence $(Y_i)_{i \geq n}$, define $\tau_0 = \inf\{i \in \mathbb{N} : U_{0,i} = \mathbf{0}\}$, and recall the definition of the mixing coefficients $\alpha(n)$ from (1.6). On one hand, using inequality (1.8), we obtain for any $n \in \mathbb{N}$, $A \in \sigma(Y_0)$ and $B \in \sigma(Y_i : i > n)$,

$$\mathbb{P}(Y_0 \in A, \mathfrak{Y}_n \in B) \geq \mathbb{P}(Y_0 \in A, \tau_0 \leq n/2, \mathfrak{Y}_n \in B)$$
$$= E_P(P_\omega(Y_0 \in A, \tau_0 \leq n/2) P_\omega(\mathfrak{Y}_n \in B))$$
$$\geq \mathbb{P}(Y_0 \in A, \tau_0 \leq n/2) \cdot \mathbb{P}(\mathfrak{Y}_n \in B) - 4\alpha([n/2])$$
$$\geq \mathbb{P}(Y_0 \in A) \cdot \mathbb{P}(\mathfrak{Y}_n \in B) - \mathbb{P}(\tau_0 > n/2) - 4\alpha([n/2]),$$

where $[n/2]$ stands for the integer part of $n/2$. On the other hand,

$$\mathbb{P}(Y_0 \in A, \mathfrak{Y}_n \in B) \leq \mathbb{P}(Y_0 \in A, \tau_0 \leq n/2, \mathfrak{Y}_n \in B) + \mathbb{P}(\tau_0 > n/2)$$
$$= E_P(P_\omega(Y_0 \in A, \tau_0 \leq n/2) P_\omega(\mathfrak{Y}_n \in B)) + \mathbb{P}(\tau_0 > n/2)$$
$$\leq \mathbb{P}(Y_0 \in A, \tau_0 \leq n/2) \cdot \mathbb{P}(\mathfrak{Y}_n \in B)$$
$$\quad + 4\alpha([n/2]) + \mathbb{P}(\tau_0 > n/2)$$
$$\leq \mathbb{P}(Y_0 \in A) \cdot \mathbb{P}(\mathfrak{Y}_n \in B) + 4\alpha([n/2]) + \mathbb{P}(\tau_0 > n/2).$$

Thus it remains to show that $\limsup_{n\to\infty} 1/n \log \mathbb{P}(\tau_0 > n) < 0$. Using Chebyshev's inequality in the third step and Jensen's inequality in the last one, we obtain

$$\mathbb{P}(\tau_0 > n) = \mathbb{P}\left(\bigcup_{i=1}^{d}\{U_{0,n}^{(i)} \geq 1\}\right) \leq \sum_{i=1}^{d} \mathbb{P}(U_{0,n}^{(i)} \geq 1) \leq \mathbb{E}(\|U_{0,n}\|_1)$$
$$= E_P(\|M_{n-1}M_{n-2}\cdots M_0 I_0\|_1)$$
$$\leq E_P(\|M_{n-1,\kappa}M_{n-2,\kappa}\cdots M_{0,\kappa}I_{0,\kappa}\|_1),$$

yielding the claim in virtue of assumption (B1) since the matrix norm $\|\cdot\|_1$ and $\|\cdot\|_\kappa$ are equivalent. $\square$

We next check that the second moment condition of Lemma 3.1 holds under the assumptions of Theorem 1.3.

LEMMA 3.3. *Suppose that assumption* (B1) *of Theorem* 1.3 *is satisfied. Then* $\mathbb{E}(\|Y_0\|_\kappa^\kappa) < \infty$.

PROOF. Recall the definition of the norm $\|x\|_{\omega,\kappa}$, $x \in \mathbb{R}$, from (1.4). For $i = 1,\ldots,d$ we obtain, by using Minkowski's inequality,

$$\mathbb{E}((Y_0^{(i)})^\kappa) = E_P(\|Y_0^{(i)}\|_{\omega,\kappa}^\kappa) = E_P\left(\left\|\sum_{n=1}^{\infty} U_{0,n}^{(i)}\right\|_{\omega,\kappa}^\kappa\right)$$
$$\leq E_P\left(\left[\sum_{n=1}^{\infty} \|U_{0,n}^{(i)}\|_{\omega,\kappa}\right]^\kappa\right) \leq \left(\sum_{n=1}^{\infty} [E_P(\|U_{0,n}^{(i)}\|_{\omega,\kappa}^\kappa)]^{1/\kappa}\right)^\kappa.$$

Conditioning on $U_{0,n-1}$ and using Minkowski's inequality again we obtain:

$$\|U_{0,n}^{(i)}\|_{\omega,\kappa} = \left\|\sum_{j=1}^{d}\sum_{m=1}^{U_{0,n-1}^{(j)}} L_{n-1,m}^{(j,i)}\right\|_{\omega,\kappa} \leq \sum_{j=1}^{d} M_{n-1,\kappa}^{(j,i)} \|U_{0,n-1}^{(j)}\|_{\omega,\kappa},$$

where the random vectors $L_{n,m}^{(i)}$ are the same as in the branching equation (1.3). It follows by using induction that

$$E_P(\|U_{0,n}^{(i)}\|_{\omega,\kappa}^\kappa) \leq E_P(\|M_{n-1,\kappa}\cdots M_{0,\kappa}I_0\|_\kappa^\kappa),$$

and thus the claim of the lemma follows from assumption (B1). $\square$

Let $T_n = \sum_{i=0}^{n-1} Y_i$ and recall the vector $\rho$ defined in (1.5).

LEMMA 3.4. *Let Assumption* 1.1 *and conditions* (B1) *and* (B2) *of Theorem* 1.3 *hold. Then the sequence* $n^{-1/2}(T_n - n \cdot \rho)$ *converges in law to a Gaussian random vector* $T_\infty$ *with zero mean. Moreover, if condition* (B3$^{(i)}$) *holds for some* $i \in \{1,\ldots,d\}$ *then* $\mathbb{E}([T_\infty^{(i)}]^2) > 0$ *for that* $i$.

PROOF. The vectors $Y_n$ satisfy conditions of Lemma 3.1 in view of Lemmas 3.2 and 3.3. Thus we only need to check that the $i$th component of the limit random vector $T_\infty$ is nondegenerate if condition (B3$^{(i)}$) holds.

Let $R_n(\omega) = E_\omega(Y_n)$. By (3.3) and using the fact that $Y_0$ and $Y_n$ are independent in a fixed environment, we obtain for any $i = 1, \ldots, d$,

$$\mathbb{E}([T_\infty^{(i)}]^2) = \mathbb{E}([Y_0^{(i)} - \rho^{(i)}]^2) + 2\sum_{n=1}^\infty \mathbb{E}[(Y_0^{(i)} - \rho^{(i)})(Y_n^{(i)} - \rho^{(i)})]$$

$$= \mathbb{E}([Y_0^{(i)}]^2) - E_P([R_0^{(i)}]^2) + E_P([R_0^{(i)}]^2) - (\rho^{(i)})^2$$

$$+ 2\sum_{n=1}^\infty E_P[(R_0^{(i)} - \rho^{(i)})(R_n^{(i)} - \rho^{(i)})]$$

$$:= \mathbb{E}([Y_0^{(i)}]^2) - E_P([R_0^{(i)}]^2) + \sigma_i.$$

By Lemma 3.1 the following series converges absolutely:

$$\sum_{n=1}^\infty \mathbb{E}[(Y_0^{(i)} - \alpha_3^{(i)})(Y_n^{(i)} - \rho^{(i)})] = \sum_{n=1}^\infty E_P[(R_0^{(i)} - \rho^{(i)})(R_n^{(i)} - \rho^{(i)})].$$

Moreover, it follows from Minkowski's inequality that

$$E_P([R_0^{(i)}]^2) = E_P\left(\left[\sum_{n=1}^\infty (M_{n-1}\cdots M_0 I_0)^{(i)}\right]^2\right)$$

$$\leq \left[\sum_{n=1}^\infty \sqrt{E_P([(M_{n-1}\cdots M_0 I_0)^{(i)}]^2)}\right]^2$$

and hence $E_P([R_0^{(i)}]^2) < \infty$ in view of (B1). Therefore (see the last few lines in [5], page 198), $\sigma_i = \lim_{n\to\infty} \frac{1}{n} E_P([\sum_{j=0}^{n-1} R_j^{(i)} - n \cdot \rho^{(i)}]^2) \geq 0$. Furthermore, if condition (B3$^{(i)}$) holds, $Y_0^{(i)}$ is a nondegenerate random variable under $P_\omega$ for $P$-a.e. environment $\omega$. It follows that

$$\mathbb{E}([Y_0^{(i)}]^2) - E_P([R_0^{(i)}]^2) = E_P[E_\omega([Y_0^{(i)}]^2) - [E_\omega(Y_0^{(i)})]^2] > 0,$$

completing the proof of the lemma. $\square$

To complete the proof of Theorem 1.3 observe that

$$T_n - S_n = \sum_{i=0}^{n-1} \sum_{j=n}^\infty U_{i,j} \leq \sum_{i=-\infty}^{n-1} \sum_{j=n}^\infty U_{i,j}$$

and hence, in virtue of assumption (B1), $\mathbb{E}(T_n - S_n)$ is bounded by a vector in $\mathbb{R}^d$: $\mathbb{E}(\sum_{i=-\infty}^{n-1} \sum_{j=n}^\infty U_{i,j}) = \mathbb{E}(\sum_{i=-\infty}^{-1} \sum_{j=0}^\infty U_{i,j}) \in \mathbb{R}^d$. It follows that $(T_n - S_n)/\sqrt{n}$ converges to $\mathbf{0}$ in probability, yielding the claim of Theorem 1.3.

**4. Recurrence behavior of the branching process.** This section is devoted to the proof of Theorem 1.5 and Theorem 1.6. The proof of the former is included in Section 4.1 while that of the latter is contained in Section 4.2.

4.1. *Proof of Theorem* 1.5. First we make explicit the obvious link between the structure of the limiting distribution and the recurrence properties of the sequence $(Z_n)_{n \in \mathbb{Z}_+}$. For any environment $\omega$, let

$$\pi_\omega(v) = P_\omega(\widetilde{Z}_0 = v), \qquad v \in \mathbb{Z}_+^d, \tag{4.1}$$

and let $\pi(v)$ denote $E_P(\pi_\omega(v)) = \mathbb{P}(\widetilde{Z}_0 = v)$. According to Theorem 2.1, $\pi_\omega$ is a probability distribution on $\mathbb{Z}_+^d$ for $P$-almost every environment $\omega$. In the following lemma we use the notation $Z = (Z_n)_{n \geq 0}$ an $\widetilde{Z} = (\widetilde{Z}_n)_{n \geq 0}$.

LEMMA 4.1. *Let Assumption* 1.1 *hold. Then:*

(a) $\mathbb{P}(Z_n = v \text{ i.o.}) \in \{0, 1\}$ *for every* $v \in \mathbb{Z}_+^d$. *Moreover,* $\mathbb{P}(Z_n = v \text{ i.o.}) = 1$ *if and only if* $\pi(v) > 0$.

(b) *If* $\pi(\mathbf{0}) > 0$ *then* $\mathbb{P}(\widetilde{Z} \in \cdot | \widetilde{Z}_0 = \mathbf{0}) = \int_\Omega P_\omega(Z \in \cdot) \widetilde{P}(d\omega)$, *where the probability measure* $\widetilde{P}$ *on the space of environments* $(\Omega, \mathcal{F})$ *is defined by* $\frac{d\widetilde{P}}{dP} = \frac{\pi_\omega(\mathbf{0})}{\pi(\mathbf{0})}$. *In particular, the measures* $\mathbb{P}(Z \in \cdot)$ *and* $\mathbb{P}(\widetilde{Z} \in \cdot | \widetilde{Z}_0 = \mathbf{0})$ *are equivalent if and only if* $P(\pi_\omega(\mathbf{0}) > 0) = 1$.

PROOF. (a) This is a direct consequence of Lemma 2.2. By part (c) of the lemma, $\mathbb{P}(Z_n = v \text{ i.o.}) = \mathbb{P}(\widetilde{Z}_n = v \text{ i.o.})$, while part (b) implies that $\mathbb{P}(\widetilde{Z}_n = v \text{ i.o.}) \in \{0, 1\}$ and furthermore, $\mathbb{P}(\widetilde{Z}_n = v \text{ i.o.}) = 1$ if and only if $\mathbb{P}(\widetilde{Z}_0 = v) > 0$ (the latter follows from Poincaré recurrence theorem in one direction and from the Borel–Cantelli lemma in the other one).

(b) For any measurable set $A \subseteq (\mathbb{Z}_+^d)^{\mathbb{Z}_+}$ we have

$$\frac{\mathbb{P}(\widetilde{Z} \in A \cap \widetilde{Z}_0 = \mathbf{0})}{\mathbb{P}(\widetilde{Z}_0 = \mathbf{0})} = \frac{1}{\pi(\mathbf{0})} E_P[P_\omega(\widetilde{Z} \in A | \widetilde{Z}_0 = \mathbf{0}) P_\omega(\widetilde{Z}_0 = \mathbf{0})]$$

$$= E_P\left[P_\omega(Z \in A) \frac{\pi_\omega(\mathbf{0})}{\pi(\mathbf{0})}\right]. \qquad \square$$

It should be noted that $P(\pi_\omega(\mathbf{0}) > 0) \in (0, 1)$ in a fairly common situation. Consider, for example, an environment which is a Markov chain in the state space $\{\alpha, \beta\}$ and a process $Z_n$ such that $\mathbb{P}(Z_1 = \mathbf{0} | \omega_0 = \alpha) = 0$ while $\mathbb{P}(Z_1 = \mathbf{0} | \omega_0 = \beta) = 1$. However, using the fact that $\pi_{\theta^n \omega}(\mathbf{0}) \geq \pi_\omega(v) P_\omega(\widetilde{Z}_n = \mathbf{0} | \widetilde{Z}_0 = v)$ for every $v \in \mathbb{Z}_+^d$, the following simple criterion was in fact proved by Key in [14] (see the proof of the second part of [14], Theorem 3.3).

LEMMA 4.2 ([14]). *Let Assumption* 1.1 *hold and suppose in addition that for some $\omega \in \Omega$ there exists $n \in \mathbb{N}$ such that $P_\omega(\widetilde{Z}_n = \mathbf{0} | \widetilde{Z}_0 = \mathbf{1}) > 0$, where $\mathbf{1} = (1, \ldots, 1) \in \mathbb{R}^d$. Then, $\pi_{\theta^n \omega}(\mathbf{0}) > 0$ with probability one.*

We note in passing that if the environment is an i.i.d. sequence and the law of the immigration $q_n$ is independent of the branching law $p_{n-1}$ for all $n \in \mathbb{Z}$, then, integrating the equation $\pi_{\theta \omega}(v) = \sum_{u \in \mathbb{Z}_+^d} \pi_\omega(u) P_\omega(\widetilde{Z}_1 = v | \widetilde{Z}_0 = u)$, we obtain that $\pi$ is a solution of the eigen-vector problem $\pi = \pi A$ for the stochastic matrix $A_{u,v} = \mathbb{P}(\widetilde{Z}_1 = v | \widetilde{Z}_0 = u)$.

Part (a) of Lemma 4.1 combined with Lemma 4.2 imply that under the conditions of the theorem, $\mathbb{P}(Z_n = \mathbf{0} \text{ i.o.}) = 1$. The rest of the subsection is devoted to the proof that $\mathbb{P}(\lim_{n \to \infty} \nu_n / n = 1/\pi(\mathbf{0})) = 1$ if (1.11) holds.

Let $\widetilde{\mathbb{P}}$ be the probability measure on the product space $\Omega \times (\mathbb{Z}_+^d)^{\mathbb{Z}_+}$ defined by setting $\widetilde{\mathbb{P}}(\omega \in A, (Z_n)_{n \in \mathbb{Z}_+} \in B) = \mathbb{P}(\omega \in A, (\widetilde{Z}_n)_{n \in \mathbb{Z}_+} \in B | \widetilde{Z}_0 = \mathbf{0})$ for measurable sets $A \subseteq \Omega$ and $B \subseteq (\mathbb{Z}_+^d)^{\mathbb{Z}_+}$. The expectation with respect to $\widetilde{\mathbb{P}}$ is denoted by $\widetilde{\mathbb{E}}$. For $n \in \mathbb{N}$ let $\mu_n = \nu_n - \nu_{n-1}$ and consider the triples $x_n = (\theta^{\nu_n} \omega, \mu_n, Z_{[\nu_{n-1}+1, \nu_n]})$, where the shift $\theta$ is defined in (1.13) and the segments $Z_{[i,j]}$ in (1.12). The sequence $(x_n)_{n \in \mathbb{N}}$ is a Markov chain on the state space $\mathbb{X} = \Omega \times \mathbb{N} \times \Upsilon_d$, where $\Upsilon_d = \bigcup_{n \in \mathbb{Z}_+} (\mathbb{Z}_+^d \setminus \{\mathbf{0}\})^n \times \{\mathbf{0}\}$. We equip $\Upsilon_d$ with the discrete topology and denote by $\Xi$ the product $\sigma$-algebra induced on $\mathbb{X}$. Transition kernel $K$ of $x_n$ under $\widetilde{\mathbb{P}}$ is given by

$$K(\omega, i, \mathbf{x}; A, j, \mathbf{y}) = P_\omega(Z_{[1,j]} = \mathbf{y}) \mathbf{I}_{\theta^{-j} A}(\omega),$$

where $A \in \mathcal{F}, i, j \in \mathbb{N}, \mathbf{x}, \mathbf{y} \in \Upsilon_d$, and $\omega \in \Omega$. The initial law of $x_n$ is given by

$$\mathbb{Q}(A, i, \mathbf{x}) := \widetilde{\mathbb{P}}(\omega \in A, \nu_1 = i, Z_{[1, \nu_1]} = \mathbf{x})$$

$$= \int_\Omega P_\omega(Z_{[1,i]} = \mathbf{x}) \mathbf{I}_{\theta^{-i} A}(\omega) \widetilde{P}(d\omega),$$

where the measure $\widetilde{P}$ is introduced in part (b) of the statement of Lemma 4.1.

It follows from Lemmas 4.1 and 4.2 that if (1.11) holds, then $\widetilde{P}$ is equivalent to $\mathbb{P}$ and hence $\mathbb{Q}$ is equivalent to $\mathbb{P}$. Moreover, applying Kac's recurrence theorem ([10], page 348) to the stationary sequence $\widetilde{Z}_n$, we infer that $\widetilde{\mathbb{E}}(\mu_1) = 1/\pi(\mathbf{0})$. Therefore, the second part of Theorem 1.5 is implied by the following proposition.

PROPOSITION 4.3. *Let Assumption* 1.1 *and condition* (1.11) *hold. Then the Markov chain $(x_n)_{n \in \mathbb{N}}$ on $(\mathbb{X}, \Xi)$ with initial distribution $\mathbb{Q}$ and transition kernel $K$ is stationary and ergodic.*

PROOF. The conclusion that the measure $\mathbb{Q}(\cdot) = \mathbb{P}(\cdot | \widetilde{Z}_0 = \mathbf{0})$ is preserved under the action of the kernel $K$ is standard (see, e.g., [15], Chapter II or

[16], Section 2.3) and follows from the fact that $\mathbb{P}$ is preserved under the shift on the sequence $y_n = (\theta^n \omega, \widetilde{Z}_n)$. Specifically, for integers $i, j$ such that $i \leq j$ let $\widetilde{Z}_{[i,j]} = (\widetilde{Z}_i, \widetilde{Z}_{i+1}, \ldots, \widetilde{Z}_j)$ and observe that for any $i \in \mathbb{N}, \mathbf{x} \in \Upsilon_d, A \in \mathcal{F}$,

$$\mathbb{Q}(A, i, \mathbf{x}) = \frac{1}{\pi(\mathbf{0})} \mathbb{P}(\widetilde{Z}_0 = 0, \widetilde{Z}_{[1,i]} = \mathbf{x}, \theta^i \omega \in A)$$

$$= \frac{1}{\pi(\mathbf{0})} \mathbb{P}(\widetilde{Z}_{-i} = 0, \widetilde{Z}_{[-i+1,0]} = \mathbf{x}, \omega \in A).$$

Therefore, using the notation $K\mathbb{Q}(A, j, \mathbf{y}) := \int_{\mathbb{X}} K(\omega; i, \mathbf{x}; A, j, \mathbf{y}) \mathbb{Q}(d\omega, i, \mathbf{x})$,

$$K\mathbb{Q}(A, j, \mathbf{y})$$

$$= \sum_{i \in \mathbb{N}} \sum_{\mathbf{x} \in \Upsilon_d} \int_{\Omega} P_\omega(\widetilde{Z}_{-i} = 0, \widetilde{Z}_{[-i+1,0]} = \mathbf{x}) P_\omega(Z_{[1,j]} = \mathbf{y}) \mathbf{I}_{\theta^{-j}A}(\omega) \frac{dP(\omega)}{\pi(\mathbf{0})}$$

$$= \int_{\Omega} P_\omega(\widetilde{Z}_0 = 0, \widetilde{Z}_{[1,j]} = \mathbf{y}) \mathbf{I}_{\theta^{-j}A}(\omega) \frac{dP(\omega)}{\pi(\mathbf{0})} = \mathbb{Q}(A, j, \mathbf{y}).$$

Thus $\mathbb{Q}$ is an invariant measure for the transition kernel $K$.

Our next aim is to show that the sequence $(x_n)_{n \in \mathbb{N}}$ is ergodic under $\widetilde{\mathbb{P}}$. The proof is a variation of the "environment viewed from the particle" argument adapted from the theory of random motion in random media (see, e.g., [4], Lecture 1 and references therein). Let $A \in \Xi^{\otimes \mathbb{N}}$ be an invariant set, that is, $\Theta^{-1}A = A$, $\widetilde{\mathbb{P}}$-a.s., where $\Theta$ is the usual shift on the space $\mathbb{X}^{\mathbb{N}}$: $(\Theta\mathcal{X})_n = x_{n+1}$ for $\mathcal{X} = (x_n)_{n \in \mathbb{N}}$. It suffices to show that $\widetilde{\mathbb{P}}(A) \in \{0, 1\}$.

For $x \in \mathbb{X}$ let $\widetilde{\mathbb{P}}_x$ be the law of the Markov chain $(x_n)_{n \in \mathbb{N}}$ with initial state $x_1 = x$. Set $h(x) = \widetilde{\mathbb{P}}_x(A)$ and note that the sequence $(h(x_n))_{n \in \mathbb{N}}$ form a martingale in its canonical filtration. This follows from the following representation of $h(x_n)$:

$$\widetilde{\mathbb{E}}(\mathbf{I}_A(\mathcal{X})|x_0, x_1, \ldots, x_n) = \widetilde{\mathbb{E}}(\mathbf{I}_{\Theta^{-n}A}(\mathcal{X})|x_0, x_1, \ldots, x_n)$$
(4.2)
$$= h(x_n), \quad \widetilde{\mathbb{P}}\text{-a.s.,}$$

where the first equality is due to the invariance property of the set $A$ while the second one follows from the Markov property. The representation (4.2) yields, by Lévy 0–1 law (cf. [10], page 263), that $h(x_n)$ converges to $\mathbf{I}_A(\mathcal{X})$ as $n$ approaches to infinity. This in turn implies that there exists a set $B \in \Xi$ such that $h(x) = \mathbf{I}_B(x)$, $\mathbb{Q}$-a.s. (see the proof of (1.17) in [4]).

It remains to show that $\mathbb{Q}(B) \in \{0, 1\}$. By the martingale property of the sequence $h(x_n)$, $\mathbf{I}_B(x) = K\mathbf{I}_B(x)$, $\mathbb{Q}$-a.s. Since the transition kernel $K(\omega, i, \mathbf{x}; \cdot)$ does not depend of $i$ and $\mathbf{x}$, the function $h$ should not depend on these two variables as well, that is, modulo a $\mathbb{Q}$-null set, the set $B$ has the form

$$B = \{(\omega, i, \mathbf{x}) \in \mathbb{X} : \omega \in C\} = C \times \mathbb{N} \times \Upsilon_d$$

for some $C \in \mathcal{F}$.

Let $\widetilde{K}$ be the transition kernel of the Markov chain $\widetilde{\omega}_n = \theta^{\nu_n}\omega$, that is, for $i \in \mathbb{N}$ and $\mathbf{x} \in \Upsilon_d$,

$$\widetilde{K}(\omega, A) = \sum_{j \in \mathbb{N}} \sum_{\mathbf{y} \in \Upsilon_d} K(\omega, i, \mathbf{x}; A, j, \mathbf{y}) = \sum_{j \in \mathbb{N}} P_\omega(\mu_1 = j)\mathbf{I}_{\theta^{-j}A}(\omega).$$

Then

$$\widetilde{K}\mathbf{I}_C(\omega) = \sum_{j \in \mathbb{N}} \sum_{\mathbf{y} \in \Upsilon_d} K(\omega, i, \mathbf{x}; \theta^j\omega, j, \mathbf{y})\mathbf{I}_B(\theta^j\omega, j, \mathbf{y})$$

$$= \mathbf{I}_B(\omega, i, \mathbf{x}) = \mathbf{I}_C(\omega), \qquad P\text{-a.s.},$$

where the second equality follows from the fact that $K\mathbf{I}_B = \mathbf{I}_B$, $\mathbb{Q}$-a.s. Thus, we have proved that

$$\mathbf{I}_C(\omega) = \sum_{j \in \mathbb{N}} P_\omega(\mu_1 = j)\mathbf{I}_C(\theta^j\omega), \qquad P\text{-a.s.}$$

Since $P$-a.s. the random variable $\mathbf{I}_C$ takes only two values either 0 or 1, $P_\omega(\mu_1 = 1) > 0$ and $\sum_{j \in \mathbb{N}} P_\omega(\mu_1 = j) = 1$, we obtain that

$$\mathbf{I}_C(\omega) = \mathbf{I}_C(\theta\omega), \qquad P\text{-a.s.}$$

That is, $\mathbf{I}_C(\omega)$ is an invariant function of the ergodic sequence $(\omega_n)_{n \in \mathbb{Z}}$. From this it follows that $P(C) \in \{0, 1\}$ and hence $\mathbb{Q}(B) \in \{0, 1\}$. This completes the proof of the proposition. $\square$

4.2. *Proof of Theorem* 1.6. First, we construct in an enlarged probability space a representation of the law of $Z_n$ as a convex combination of an "external" multitype Galton–Watson process independent of the environment and a modified MBPIRE. The law $\mathbb{P}_e$ of the external process is defined as a probability measure on $(\mathcal{T}, \mathcal{G})$ such that $X_n \in \mathbb{Z}_+^d$ are i.i.d. with distribution $e^{(d+1)}$, $L_{n,m}^{(i)} \in \mathbb{Z}_+^d$ are i.i.d. independent of $(X_n)_{n \in \mathbb{Z}}$, and the law of $L_{n,m}^{(i)}$ is $e^i$ for all $n, m$ [the sequence $(Z_n)_{n \in \mathbb{Z}_+}$ is defined as before by equation (1.3) with initial condition $Z_0 = \mathbf{0}$]. Roughly speaking, we attach to each particle a Bernoulli random variable *independent of the environment* in such a way that if the variable takes value 1 then the particle produces progeny according to the law $\mathbb{P}_e$. Thus, *only* if the random variable associated with the particle takes value 0, the number of its progeny depends on the environment.

We next show that for arbitrary large $l_0 \in \mathbb{N}$ there is a positive probability (that does not depend on $i$) that the branching process does not use the environment in time between $\nu_i$ and $\nu_i + l_0$ and, moreover, $\nu_{i+l_0} = \nu_i + l_0$. Due to the assumed mixing property of the environment, the branching process starts after occurrence of this event in the environment $\theta^{\nu_i+l_0}\omega$ which is distributed "almost" according to $P$. The claim of Theorem 1.6 is derived then from this observation. The construction used in the proof is similar to that of $L$-safeguards for a RWRE introduced in [17], Section 3.1 and [7].

Let $\xi = (\xi_n)_{n \in \mathbb{Z}_+}$ and $\lambda = (\lambda_{n,m}^{(i)})_{n \in \mathbb{Z}_+, m \in \mathbb{N}, i=1,\ldots,d}$ be two collections of i.i.d. Bernoulli random variables on a probability space $(E, \mathcal{E}, Q_\varepsilon)$ such that $\xi$ and $\lambda$ are independent of each other under $Q_\varepsilon$ and

$$Q_\varepsilon(\xi_n = 1) = Q_\varepsilon(\lambda_{n.m}^{(i)} = 1) = \varepsilon,$$
$$Q_\varepsilon(\xi_n = 0) = Q(\lambda_{n.m}^{(i)} = 0) = 1 - \varepsilon,$$

where $\varepsilon$ is introduced in the conditions of the theorem.

Further, let $\widehat{P}$ denote the product measure on $(\Omega \times E, \mathcal{F} \times \mathcal{E})$ whose marginal on $(\Omega, \mathcal{F})$ is $P$ and that on $(E, \mathcal{E})$ is $Q_\varepsilon$. The triple $(\omega, \xi, \lambda)$ serves as a new environment for the branching process $(\widehat{Z}_n)_{n \in \mathbb{Z}_+}$ defined as follows. We assume that for any realization $(\omega, \xi, \lambda)$ of the environment, a law $P_{\omega, \xi, \lambda}$ on the underlying probability space $(\mathcal{T}, \mathcal{G})$ is defined in such a way that:

- Similarly to (1.3), $(\widehat{Z}_n)_{n \in \mathbb{Z}_+}$ is a nonhomogeneous Markov chain that satisfies initial condition $\widehat{Z}_0 = \mathbf{0}$ and branching equation

$$\widehat{Z}_{n+1} = \sum_{i=1}^{d} \sum_{m=1}^{\widehat{Z}_n^{(i)} + \widehat{X}_n^{(i)}} \widehat{L}_{n,m}^{(i)} \qquad \text{for } n \geq 0,$$

where
- $\widehat{X}_n = \mathbf{I}_{\{\xi_n = 1\}} X_{n,1} + \mathbf{I}_{\{\xi_n = 0\}} X_{n,2}$,
- $\widehat{L}_{n,m}^{(i)} = \mathbf{I}_{\{\lambda_{n,m}^{(i)} = 1\}} L_{n,m,1}^{(i)} + \mathbf{I}_{\{\lambda_{n,m}^{(i)} = 0\}} L_{n,m,2}^{(i)}$,
- each one of the four sequences $(X_{n,1})_{n \in \mathbb{Z}_+}$, $(X_{n,2})_{n \in \mathbb{Z}_+}$, $(L_{n,m,1}^{(i)})_{n \in \mathbb{Z}, m \in \mathbb{Z}_+, i=1,\ldots,d}$ and $(L_{n,m,2}^{(i)})_{n \in \mathbb{Z}, m \in \mathbb{Z}_+, i=1,\ldots,d}$ consists of independent random variables, the sequences are mutually independent under $P_{\omega, \xi, \lambda}$, and

$$\begin{cases} P_{\omega,\xi,\lambda}(L_{n,m,1}^{(i)} = v) = e_n^{(i)}(v), \\ P_{\omega,\xi,\lambda}(L_{n,m,2}^{(i)} = v) = \dfrac{p_{\omega,n}^{(i)}(v) - \varepsilon e_n^{(i)}(v)}{1 - \varepsilon}, \\ P_{\omega,\xi,\lambda}(X_{n,1} = v) = e_n^{(d+1)}(v), \\ P_{\omega,\xi,\lambda}(X_{n,2} = v) = \dfrac{q_{\omega,n}(v) - \varepsilon e_n^{(d+1)}(v)}{1 - \varepsilon}, \end{cases}$$

where the constant $\varepsilon$ and the measures $e^{(i)}$ are defined in the statement of the theorem.

Define the annealed probability measure

(4.3) $$\widehat{\mathbb{P}} = \widehat{P} \otimes P_{\omega,\xi,\lambda} = P \otimes Q_\varepsilon \otimes P_{\omega,\xi,\lambda}$$

on the measurable space $(\Omega \times E \times \mathcal{T}, \mathcal{F} \otimes \mathcal{E} \otimes \mathcal{G})$. It is easy to see that the law of the process $(\widehat{Z}_n)_{n \in \mathbb{Z}_+}$ under $\widehat{\mathbb{P}}$ coincides with that of $(Z_n)_{n \in \mathbb{Z}_+}$ under $\mathbb{P}$, while its law under $Q_\varepsilon \otimes P_{\omega,\xi,\lambda}$ is the same as that of $(Z_n)_{n \in \mathbb{Z}_+}$ under $P_\omega$. Therefore, part (a) of Theorem 1.6 is contained in the following lemma. Part (b) of the theorem follows from the lemma in virtue of Lemma 4.1(b) and Proposition 4.3.

LEMMA 4.4. *Let the conditions of Theorem 1.6 hold. Then*
$$\sup_{m \in \mathbb{N}} \sup_{A \in \widehat{\mathcal{G}}^{\nu_{n+m}}} \sup_{B \in \widehat{\mathcal{G}}_{\nu_m}, \mathbb{P}(B) > 0} \{\mathbb{P}(A|B) - \mathbb{P}(A)\} \to_{n \to \infty} 0,$$
*where we denote* $\widehat{\mathcal{G}}^n = \sigma(\widehat{Z}_i : i \geq n)$ *and* $\widehat{\mathcal{G}}_n = \sigma(\widehat{Z}_i, \omega_i : i < n)$.

PROOF. With a slight abuse of notation let
$$\nu_0 = 0 \quad \text{and} \quad \nu_n = \inf\{i > \nu_{n-1} : \widehat{Z}_i = \mathbf{0}\}.$$
Fix any $l_0 \in \mathbb{N}$. By the conditions of the theorem, $\mathbb{P}_e(\nu_{l_0} = l_0) > 0$. This implies that for some integer $R$,

(4.4) $$\mathbb{P}_e\left(\nu_{l_0} = l_0 \text{ and } \max_{0 \leq t \leq l_0 - 1} \|X_{t,1}\|_1 \leq R\right) > 0.$$

Define the sequence of independent under $\widehat{\mathbb{P}}$ events $(H_i)_{i \in \mathbb{Z}_+}$ by setting
$$H_i = \Big\{\xi_{\nu_i + t} = 1; \lambda^{(i)}_{\nu_i + t, m} = 1 \text{ for } t = 0, \ldots, l_0 - 1, i = 1, \ldots, d,$$
$$m = 1, \ldots, R; \nu_{i+l_0} = \nu_i + l_0; \max_{\nu_i \leq j \leq \nu_{i+l_0} - 1} \|X_{j,1}\|_1 \leq R\Big\}.$$

Note that on the event $H_i$, the branching process does not use the environment during $l_0 = \nu_{i+l_0} - \nu_i$ steps after $\nu_i$. It follows from (4.4) that the probability of $H_i$, $\widehat{\mathbb{P}}(H_i)$, is positive and does not depend on $i$.

For any $n \in \mathbb{N}$ define $m_n(l_0) = \min_{0 \leq i \leq n - l_0}\{i : H_i \text{ occurs}\}$, with the convention that the minimum over an empty set is infinity. For any $i < n - l_0$ and $t \in \mathbb{N}$, let $D_{i,t} = \{m_n(l_0) = i, \nu_i = t\}$ and
$$D_i = \bigcup_{t \in \mathbb{N}} D_{i,t} = \{m_n(l) = i\}$$
and
$$D = \bigcup_{i=0}^{n-l_0} D_i = \{m_n(l_0) < \infty\}.$$
The events $D_i$ are disjoint and can be represented in the form
$$D_i = H_0^c \cap \cdots \cap H_{i-1}^c \cap H_i.$$

Since $\widehat{\mathbb{P}}(H_i) =: p > 0$, we obtain that $\widehat{\mathbb{P}}(D) = \sum_{i=0}^{n-l_0-1}(1-p)^i p \to_{n\to\infty} 0$.

Let $\Theta$ be the shift by $\nu_1$ on the sequence $\widehat{Z} = (\widehat{Z}_n)_{n\in\mathbb{Z}_+}$, that is, $(\Theta^n \widehat{Z})_i = \widehat{Z}_{\nu_n+i}$. For $t \in \mathbb{Z}$ denote $\widehat{\mathbb{P}}^t = P \otimes Q_\varepsilon \otimes P_{\theta^t\omega,\xi,\lambda}$, where $\theta$ is the shift operator defined in (1.13). Using inequality (1.14), we obtain that for any events $A \in \widehat{\mathcal{G}}^n$ and $B \in \mathcal{G}_0$:

$$\widehat{\mathbb{P}}(A|B) \leq \widehat{\mathbb{P}}(A, D|B) + \widehat{\mathbb{P}}(D^c|B) = \sum_{i=0}^{n-l_0}\sum_t \widehat{\mathbb{P}}(A, D_{i,t}|B) + \widehat{\mathbb{P}}(D^c)$$

$$= \sum_{i=0}^{n-l_0}\sum_t \widehat{\mathbb{P}}(A|D_{i,t}, B)\widehat{\mathbb{P}}(D_{i,t}) + \widehat{\mathbb{P}}(D^c)$$

$$\leq \sum_{i=0}^{n-l_0}\sum_t (\widehat{\mathbb{P}}^{t+l_0}(\Theta^{-i-l_0}A) + 2\varphi(l_0))\widehat{\mathbb{P}}(D_{i,t}) + \widehat{\mathbb{P}}(D^c)$$

$$\leq \sum_{i=0}^{n-l_0} \widehat{\mathbb{P}}(\Theta^{-i-l_0}A)\widehat{\mathbb{P}}(D_i) + 2\varphi(l_0) + \widehat{\mathbb{P}}(D^c).$$

Similarly,

$$\widehat{\mathbb{P}}(A|B) \geq \widehat{\mathbb{P}}(A, D|B) = \sum_{i=0}^{n-l_0}\sum_t \widehat{\mathbb{P}}(A, D_{i,t}|B)$$

$$= \sum_{i=0}^{n-l_0}\sum_t \widehat{\mathbb{P}}(A|D_{i,t}, B)\widehat{\mathbb{P}}(D_{i,t})$$

$$\geq \sum_{i=0}^{n-l_0}\sum_x (\widehat{\mathbb{P}}^{t+l_0}(\Theta^{-i-l_0}A) - 2\varphi(l_0))\widehat{\mathbb{P}}(D_{i,t})$$

$$\geq \sum_{i=0}^{n-l_0} \widehat{\mathbb{P}}(\Theta^{-i-l_0}A)\widehat{\mathbb{P}}(D_i) - 2\varphi(l_0).$$

It follows that $|\widehat{\mathbb{P}}(A|B) - \widehat{\mathbb{P}}(A)| \leq 4(\varphi(l_0) + \widehat{\mathbb{P}}(D^c))$. That is, for $n$ large enough we have $|\widehat{\mathbb{P}}(A|B) - \widehat{\mathbb{P}}(A)| \leq 8\varphi(l_0)$, proving the lemma since $l_0$ is an arbitrary integer. $\square$

**Acknowledgments.** I am grateful to Professor A.-S. Sznitman for helpful and stimulating discussions on the subject of this paper. I would like to thank the FIM at ETH Zurich for hospitality and financial support during a visit in which part of this work was carried out.

DEPARTMENT OF MATHEMATICS
UNIVERSITY OF BRITISH COLUMBIA
121-1984 MATHEMATICS ROAD
VANCOUVER, BRITISH COLUMBIA
CANADA V6T 1Z2
E-MAIL: roiterst@math.ubc.ca
URL: www.math.ubc.ca/~roiterst